\documentclass[12pt]{amsart}

\usepackage{amsmath,amsthm,amscd,euscript}
\def\s{\sigma}

\def\r{\mathbb R}

\title[Multiboundary singularities $B_n^l$
and Bernoulli--Euler numbers]{Combinatorics of multiboundary
singularities $B_n^l$ and Bernoulli--Euler numbers.}
\author{Karpenkov~O.}
\thanks{This paper was partially supported by RBRF-00-15-96084
and RBRF-01-01-00803 projects}
\date{}

\begin{document}
\maketitle

\begin{center}
{\bf Introduction.}
\end{center}

In [1,2] V.~I.~Arnold in particular established the connection
between the components of the space of very nice M-morsifications
for boundary singularities $B_n$ and the combinatorics of
corresponding Springer cones; the number of the components equals
Bernoulli--Euler number $K_n$.

In this note we regard the generalization of the boundary
singularities $B_n$ of the functions on the real line to the case
where the boundary is a finite number of ($l$) points. These
singularities $B_n^l$ could also arise in higher dimensional
case, when the boundary is an immersed hypersurface.

We obtain some recurrence relation on the numbers of connected
components of very nice M-morsification spaces with different
values of $n$ and $l$.

The author is grateful to V.~I.~Arnold, V.~M.~Zakalukin and
S.~K.~Lando for useful discussions and attention to this work.

\begin{center}
{\bf Main notions and definitions.}
\end{center}

{\bf Definition 1.} {\it A very nice M-morsification of a
multiboundary singularity $B_n^l$} is a polynomial with all its
critical points being real. All critical values and all values at
the boundary points $x=b_i$ are also different.

Note that we enumerate the boundary points, otherwise we regard
the factorization of $\r^k$ over the action of the group of
coordinate permutations. We enumerate the boundary points, since
they correspond to the different preimages, and this preimages
don't permute.

{\bf Definition 2.} The {\it M-domain} is a closed subset of the
space of polynomials
$$
x^n+\lambda _2 x^{n-2}+ \cdots +\lambda_{n-1}x
$$
consisting of polynomials which critical points are real.

The set of very nice M-morsifications is an open set in $\r
^{n-2} \times  \r ^l$. The closure of this set is $(${\it
M-domain}$) \times \r ^l$. It is subdivided into connected
components by the bifurcation diagram, containing five
hypersurfaces. Three of these hypersurfaces may occur in the case
of the ordinary boundary singularities $B_n$ (see also [2]):
\begin{description}
\item[(a)] the boundary caustic consisting of functions with a boundary
critical point;
\item[(b)] the ordinary Maxwell stratum consisting of functions with
equal critical values at different points;
\item[(c)] the boundary Maxwell stratum consisting of functions having
some value at the boundary equals some critical value (the
corresponding critical point is not at the boundary).
\end{description}
Notice that $B_n=B_n^1$. Moreover, the boundary point does not
fixed. So the definition of the very nice M-morsification of the
boundary singularity $B_n$ in the paper [2] is the special case
of Definition~1.

Two new hypersurfaces occur
in the case of the multiboundary singularities $B_n^l$, $l\ge 2$:
\begin{description}
\item[(d)] the double boundary caustic consisting of functions with
some double boundary point;
\item[(e)] the double boundary Maxwell stratum consisting of functions
with equal values at different boundary points.
\end{description}

Consider an example of the multiboundary singularity $B_3^2$
(Fig.~1) Here we regard the family of polynomials $x^3+\lambda x$
with boundary points $x=b_1$ and $x=b_2$. The M-domain is a
half-space with coordinates $(\lambda, b_1, b_2)$, $\lambda \le
0$. The hypersurfaces are marked with the following symbols: $a$
--- the boundary caustic, $c$ --- the boundary Maxwell stratum,
$d$ --- the double boundary caustic, $e$ --- the double boundary
Maxwell stratum, $f$ --- the ordinary caustic.

There is no Maxwell stratum here, since for any polynomial of the
third power with real critical values the equivalence of critical
values implies the equivalence of critical points.

\begin{figure}[tbp]
\unitlength 1.00mm \linethickness{0.4pt}
\begin{picture}(131.67,138.00)
\put(50.03,90.10){\line(0,-1){44.25}}
\put(25.28,54.53){\line(0,1){71.92}}
\put(44.98,35.59){\line(0,1){69.87}}
\put(50.10,29.96){\line(0,1){70.13}}
\put(104.87,54.79){\line(0,-1){39.67}}
\multiput(79.79,115.19)(0.12,-0.12){209}{\line(1,0){0.12}}
\put(104.87,90.11){\line(0,-1){35.07}}
\multiput(104.87,15.12)(-0.12,0.32){113}{\line(0,1){0.32}}
\multiput(79.91,115.06)(-0.12,0.12){160}{\line(-1,0){0.12}}
\put(60.73,134.24){\line(-1,0){43.41}}
\multiput(60.73,133.89)(0.12,-0.33){56}{\line(0,-1){0.33}}
\put(10.67,80.00){\line(1,0){34.33}}
\put(131.67,80.00){\vector(1,0){0.2}}
\put(105.33,80.00){\line(1,0){26.34}}
\put(80.00,129.33){\vector(0,1){0.2}}
\put(80.00,80.00){\line(0,1){49.33}}
\put(80.00,15.00){\line(0,-1){4.33}}
\multiput(110.00,50.00)(-0.12,0.12){250}{\line(0,1){0.12}}
\put(25.00,137.67){\vector(-1,1){0.2}}
\multiput(28.33,134.33)(-0.12,0.12){28}{\line(0,1){0.12}}
\put(128.33,83.33){\makebox(0,0)[cc]{$\lambda$}}
\put(29.33,138.00){\makebox(0,0)[cc]{$b_2$}}
\put(72.00,5.00){\makebox(0,0)[cc] {Fig.1 Bifurcation diagram of
the M-domain for the singularity $B_3^2$ $(\lambda \le 0)$.}}
\put(69.33,21.33){\makebox(0,0)[cc]{$d$}}
\put(105.00,15.00){\line(-1,0){43.67}}
\put(30.33,121.67){\line(0,-1){73.00}}
\bezier{252}(25.33,78.67)(25.67,113.00)(45.00,92.33)
\bezier{240}(30.33,65.33)(50.00,46.33)(50.00,79.33)
\bezier{60}(45.00,92.00)(49.33,87.67)(50.00,79.00)
\bezier{60}(25.33,78.67)(25.67,71.33)(30.33,65.33)
\bezier{188}(65.00,39.33)(72.33,51.33)(105.00,55.00)
\bezier{176}(105.00,55.00)(81.67,52.67)(65.00,64.67)
\bezier{212}(65.00,30.67)(70.33,47.67)(105.00,55.00)
\bezier{192}(105.00,55.00)(82.33,54.67)(65.00,72.67)
\bezier{228}(25.33,126.33)(51.33,114.00)(80.00,114.67)
\bezier{152}(80.00,114.67)(53.67,110.67)(50.00,100.00)
\bezier{208}(30.33,121.67)(49.00,113.00)(80.00,114.67)
\bezier{160}(80.00,114.67)(45.67,110.67)(45.00,105.33)
\bezier{128}(50.00,88.00)(60.33,82.67)(80.00,79.67)
\bezier{28}(45.00,92.67)(50.00,89.33)(50.00,90.00)
\bezier{72}(30.33,107.67)(37.67,101.00)(45.00,96.67)
\bezier{48}(15.00,122.67)(19.67,119.00)(24.33,115.33)
\bezier{24}(25.00,112.33)(29.00,109.33)(30.33,109.00)
\bezier{80}(50.00,54.00)(51.00,62.33)(60.00,69.00)
\bezier{28}(50.00,45.33)(51.33,49.00)(53.00,51.00)
\bezier{36}(50.00,79.33)(54.33,78.33)(59.00,78.67)
\bezier{52}(15.00,114.33)(20.67,109.33)(25.33,107.00)
\bezier{72}(30.33,99.00)(37.67,93.00)(45.00,89.33)
\bezier{36}(30.33,74.00)(33.00,76.33)(38.33,77.33)
\bezier{24}(30.33,65.67)(32.33,67.00)(36.00,68.00)
\bezier{20}(47.33,62.00)(49.33,64.67)(50.00,65.67)
\bezier{20}(34.67,100.00)(39.33,98.33)(39.33,98.33)
\bezier{12}(25.33,70.33)(26.00,71.67)(27.00,71.33)
\bezier{28}(50.00,30.00)(51.00,34.00)(53.67,35.33)
\bezier{24}(45.00,35.67)(46.00,37.67)(50.00,38.33)
\bezier{64}(30.33,48.67)(35.00,45.33)(45.00,44.67)
\bezier{24}(25.33,54.67)(26.33,53.33)(30.00,52.00)
\multiput(65.33,72.33)(-0.12,0.12){426}{\line(0,1){0.12}}
\multiput(65.33,64.33)(-0.12,0.12){423}{\line(0,1){0.12}}
\multiput(65.00,39.00)(-0.12,0.12){417}{\line(0,1){0.12}}
\bezier{24}(15.00,80.67)(16.00,83.67)(18.33,85.33)
\multiput(65.33,30.33)(-0.12,0.12){420}{\line(0,1){0.12}}
\bezier{52}(15.33,89.00)(16.33,92.67)(25.33,94.67)
\multiput(61.33,15.00)(-0.12,0.32){95}{\line(0,1){0.32}}
\multiput(50.00,45.67)(-0.12,0.33){164}{\line(0,1){0.33}}
\multiput(30.33,99.67)(-0.12,0.32){109}{\line(0,1){0.32}}
\put(79.67,80.67){\line(0,1){0.00}}
\multiput(80.33,79.33)(-0.11,0.22){3}{\line(0,1){0.22}}
\put(84.67,125.67){\makebox(0,0)[cc]{$b_1$}}
\put(28.00,118.33){\makebox(0,0)[cc]{$c$}}
\put(36.00,112.00){\makebox(0,0)[cc]{$a$}}
\put(47.33,99.00){\makebox(0,0)[cc]{$a$}}
\put(56.33,91.67){\makebox(0,0)[cc]{$c$}}
\put(90.33,87.00){\makebox(0,0)[cc]{$f$}}
\put(70.33,74.00){\makebox(0,0)[cc]{$c$}}
\put(68.67,65.33){\makebox(0,0)[cc]{$a$}}
\put(68.67,50.00){\makebox(0,0)[cc]{$a$}}
\put(68.33,40.33){\makebox(0,0)[cc]{$c$}}
\put(27.67,85.67){\makebox(0,0)[cc]{$e$}}
\bezier{24}(45.00,59.67)(46.00,62.67)(48.33,64.00)
\bezier{20}(45.00,84.67)(47.67,83.00)(49.33,82.67)
\bezier{20}(25.33,78.67)(28.33,79.00)(30.33,79.00)
\bezier{20}(25.33,104.33)(28.67,102.00)(29.33,102.00)
\bezier{16}(45.00,50.67)(46.00,53.00)(47.00,54.00)
\end{picture}

\end{figure}

\begin{center}
{\bf The recurrent equation and some its corollaries}.
\end{center}

By $K_n$ we denote Bernoulli--Euler numbers (here is the beginning
of this sequence, the first number corresponds to $n=0$:
$1,1,1,2,5,16,61,\ldots$). By $K_n^l$ we denote the quantity of
connected components of the very nice M-morsification spaces for
the singularity $B_n^l$. Bernoulli--Euler numbers are the
boundary conditions for the numbers $K_n^l$, that is
$K_n^0=Z_{n-2}$ and $Z_n^1=K_{n+1}$ (See [1]).

{\bf Theorem.} The following equation on the numbers $K_n^l$
holds:
$$
K_{n-2}^{l+1}=K_n^l-nlK_n^{l-1}.
$$

To proof this theorem we need to use the $\s$-shaped functions,
as it makes in case of the boundary singularities $B_k$ (See [1]).
In the paper [1] to any polynomial with boundary point adds the
$\s$-shaped function concentrated in some neighborhood of the
boundary point. It follows that there exists a one-to-one
correspondence with the polynomials which has no boundary points
and the degree of this polynomials is greater by two. Presence of
another boundary points rather complicates the picture: the
problem loses "up-down" symmetry. Thus, we either add or
subscribe some $\s$-shaped function concentrated in the
neighborhood of some point. We regard the case where some
boundary value is greater than any critical value or conversely
smaller than any critical value. Thus, we divide the proof on the
cases of the polynomials of even and odd degree.

Finally, let us regard some corollaries of the theorem.

The exponential generating function for Bernoulli--Euler numbers
is the function $K(t)=\tan(t)+\sec(t)$. Notice that for $l=0$ and
$l=1$ we may consider $K(t)$ as an exponential generating
function, however with our notations it would be $K_0(t)=\int
K(t)dt= -\ln(\cos(t))+\ln(\tan(\frac{t}{2}+\frac{\pi}{4}))+C$ and
$K_1(t)=K'(t)=\frac{1+K^2}{2}=\frac{1+\sin(t)}{\cos^2{t}}=\frac{1}{1-\sin(t)}$
respectively.

{\bf Corollary 1.} The exponential generating functions for
$l=7,3,4$ are the following:
$$
\begin{array}{l}
K_2(t)=
\frac{3\sin(t)-t\cos(t)}{(1-\sin(t))^2};\\
K_3(t)=
\frac{7}{(1-\sin(t))^3}\Bigl(\sin(t)(3\sin(t)+7)-t\cos(t)(5+\sin(t))\Bigr);\\
K_4(t)=
\Bigl(\frac{3t^2}{1-\sin(t)}-\frac{3t\cos(t)}{(1-\sin(t))^2}(3-\sin(t))+
\frac{3(2-\sin(t))}{(1-\sin(t))^2} \Bigr)^{(4)}.\\
\end{array}
$$

Consider the exponential generating function of two variables
$K(x,y)=\sum\frac{K_n^l}{l!n!}x^ly^n$.

{\bf Corollary 2.}
$K(x,y)$ satisfies the following differential equation:
$$
K_x=(1-2x)K_{yy}-xyK_{yyy}.
$$

At last we calculate the numbers
$K_n^l$ for $n\le 0$ using the relation of the theorem.

We may regard $K_0^l$ as the number of connected components
of the space of constant polynomials with $l$ boundary points,
where any polynomial has different values.
Here all points of the real line are critical. So
the critical value equals to any boundary value.
We do not know what is the meaning of the numbers
$K_n^l$ for $n<0$.

\begin{center}
\begin{tabular}{|c|c|c|c|c|c|}
\hline
$K_n^l$ & l=1 & l=2 & l=3 & l=4 & l=5 \\
\hline
n=0     &  1  &  0  &  0  &  0  &  0  \\
\hline
n=-1    &  1  &  0  &  0  &  0  &  0  \\
\hline
n=-2    &  ?  &  1  &  0  &  0  &  0  \\
\hline
n=-3    &  ?  &  ?  &  2  &  0  &  0  \\
\hline
n=-4    &  ?  &  ?  &  ?  &  6  &  0  \\
\hline
n=-5    &  ?  &  ?  &  ?  &  ?  &  24 \\
\hline
\end{tabular}
\end{center}

Note that in any row all values are equal to zero starting
from some certain index.
Let us state this in the following corollary.

{\bf Corollary 3.}
Let $n \le -1 $, then $K_n^l=0$ for $l > -n$, and $K_n^{-n}=(-n-1)!$.
$K_0^l=0$ for $l >1$.

The proof uses the statement of the theorem and carrying out by
induction on $n$.

\begin{center}
{\bf References.}
\end{center}

[1] V.~I.~Arnold, {\it Bernoulli-Euler updown numbers associated with function
singularities, their combinatorics and arithmetics}, Duke Math. J., 1991,
63(2), 537-555.

[2] V.~I.~Arnold, {\it Springer numbers and morsification spaces}, J. Algebraic
Geom., 1992, 1(2), 197-214.
\end{document}